\documentclass[twoside,reqno,11pt]{fcaa-var} %

\usepackage{graphicx}
\usepackage{epsfig}
\usepackage{amsthm}
\usepackage{amsmath}
\usepackage{latexsym}
\usepackage{amsfonts}
\usepackage{amssymb}

\newtheoremstyle{theorem}
  {15pt}          
  {15pt}  
  {\sl}  
  {\parindent}
  {\sc}  
  {. }   
  { }    
  {}     
\theoremstyle{theorem}

\newtheoremstyle{defi}
  {15pt}          
  {15pt}  
  {\rm}  
  {\parindent}     
  {\sc}  
  {. }    
  { }    
  {}     
\theoremstyle{defi}

 \def\proofname{\indent\rm P\ r\ o\ o\ f.}
 \def\proofend{\hfill$\Box$}
 \def\qed{\hfill$\Box$} 

 \def\theequation{\arabic{section}.\arabic{equation}}

  \usepackage[colorlinks=true,citecolor=red,filecolor=red,urlcolor=magenta]{hyperref}
 
 \newtheorem{thm}{Theorem}[section]
\newtheorem{prop}{Proposition}  [section]
\newtheorem{defn}{Definition}[section]
\newtheorem{coro}{Corollary} [section]
\newtheorem{rem}{Remark} [section]
\newtheorem{exa}{Example} [section]
\newtheorem{lemme}{Lemma}[section]
 \numberwithin{equation}{section}
\def\set{[a,b]}

\def\BV{BV(\set) }
\def\Cone{\mathcal{C}^{1}(\set) }

\def\intset{\int_{\set}}

 \def\themycopyrightfootnote{\vspace*{3pt}
 \copyright \, Year\,  Diogenes  Co., Sofia
   \par  \noindent pp. xxx--xxx, DOI: ..................
   \hfill  \vspace*{-36pt}
  }

\newcommand{\R}{\mathbb{R}}

 \newcommand{\DD}{{\bf D}}
 \newcommand{\W}{\Omega}
  
\newcommand{\N} {\mathbb{N}}

\newcommand{\wspG}{W_G^{s,p}(a,b)}

\newcommand{\wsoneRL}{W^{s,1}(a,b)}
\newcommand{\wszero}{W^{s,1}_0(a,b)}

\newcommand{\wsp}{W_{RL,a+}^{s,1}(a,b)}
\newcommand{\wsm}{W_{RL,b-}^{s,1}(a,b)}
\newcommand{\lp}{L^{p}(a,b)}
\newcommand{\lone}{L^{1}(a,b)}
\newcommand{\linf}{L^{\infty}(a,b)}
\newcommand{\bv}{BV(a,b)}

 \title[Fractional Sobolev Spaces\dots]{Fractional Sobolev Spaces and   Functions  of  Bounded Variation\\ [4pt] IN ``FCAA'' JOURNAL }
\author[\normalsize M. Bergounioux, A. Leaci, G. Nardi, F. Tomarelli]{\normalsize Ma\" \i tine Bergounioux$^1$, Antonio  Leaci$^2$, Giacomo Nardi$^3$, Franco Tomarelli$^4$}

\begin{document}
\vbox to 2.5cm { \vfill }


 \bigskip \medskip

\begin{abstract}  We  investigate  the 1D Riemann-Liouville {fractional derivative }
focusing on the connections with fractional  Sobolev spaces, the space $BV$ of functions of bounded variation, whose derivatives are not functions but measures and
the space $SBV$, say the space of bounded variation functions whose derivative has no Cantor part.  We prove that $SBV$ is included in $W^{s,1} $ for every $s \in (0,1)$ while the result remains open for  $BV$. We study examples and address open questions. 

 \medskip


 \smallskip


 \end{abstract}
  \maketitle
   \vspace*{-16pt}
  
 \section{Introduction}

The aim of this work is to investigate  the  \textit{fractional derivative } concepts and make the connection between  the related  (so called fractional) Sobolev  spaces and the spaces of functions of bounded variation whose derivatives are not functions but measures.  Here, we only deal with the 1D case and  results can be  directly extended to higher dimensions via slice theorems techniques,which are not straightforward. There are two main definitions of  fractional differentiation  whose connections are not clear   (to our knowledge). The  fractional derivative in the sense of Gagliardo is not explicitely defined   (almost everywhere for example) but  through the setting of fractional Sobolev Spaces and the underlying norm (see \cite{cafarelli} for example).  It is, in  some sense, a global definition which can be easily handled via the Fourier transform in the Hilbertian case.  The second approach is based on the Riemann-Liouville fractional derivative  (in short  RL) and  may be pointwise defined. We choose to focus on the RL derivative : there are many variants of  the  fractional  derivatives/integrals  definition as
 the Grunwald-Letnikov,  Caputo,  Weyl   ones \cite{butzer00,das,Mathai}  but the  RL derivative   can be considered as a generic one.  For a complete study of these derivatives one can refer  to the book by Samko and al. \cite{samko} that contains an extensive bibliography in particular with respect to the pioneer work of Hardy-Littlewood. Moreover, in   \cite{Tarasov},   the connection  is made with metric and measure spaces, in particular the Hausdorff measure. We decided to use this derivative concept  because it seems more adapted to applications than the Gagliardo one.   The RL derivative is widely used by physicists \cite{Uchaikin_Vol1,Uchaikin_Vol2}, in  automatics, control theory and   image processing as well, especially to deal with image enhancement and texture analysis (see \cite{Siarry14b} for example); in  \cite{Malinowska12} calculus of variations problems where the cost functional involves fractional derivative are investigated.
 Nevertheless, the context  is often a discrete one and there is not much analysis (to our knowledge) in the infinite dimensional setting. In particular, the link between the classical spaces of bounded variation functions and the fractional Sobolev spaces is not clear. To our knowledge, there is no paper  that compare  the BV space and the   fractional Sobolev spaces   in  the RL sense.  Indeed,  the concept of fractional Sobolev spaces is not much  developed  for the RL derivative, though this fractional derivative concept is commonly used in engineering. One can refer  to \cite{Bourdin15,Idczak12,Idczak13} however.
 Here we consider the 1D case and
  the paper is organized as follows. Section 2 is devoted to the presentation of the two main approaches with a special focus on the Riemann-Liouville fractional derivatives: the main tools are recalled. In Section 3, we define RL-fractional Sobolev spaces  $ \wsoneRL$ and give basic properties. In the last section, we perform a comparison between these fractional Sobolev spaces, the classical  $BV$ space
 and the space $SBV$of functions whose distributional derivative is a special measure in the sense of De Giorgi (see \cite{ambrosio}). In particular we prove that
 \begin{equation*}
 SBV\subset \underset{s\in(0,1)}{\bigcap} W^{s,1}  \qquad \hbox{and}\qquad \underset{s\in(0,1)}{\bigcap} W^{s,1}\setminus BV \neq \emptyset\ .
 \end{equation*}

\section{Fractional Calculus and Fractional Sobolev Spaces}

In this section we present the  two main (different) definitions of fractional Sobolev spaces  that we can find in the literature.
We are in particular interested to the case where the differentiation order is $s \in [0,1) $ in order to study the fractional spaces between $L^1$ and $W^{1,1}$ and their relationship with $BV$.
In the sequel, we consider the 1D framework.
We recall that the space $AC(a,b)$ of absolutely continuous functions coincides with the Sobolev space $W^{1,1}(a,b)$ defined by
$$W^{1,1}(a,b):= \left \{ u \in L^{1}(a,b) ~| ~u'\in L^{1}(a,b) \, \right\}.$$
Therefore
$$W^{1,1}(a,b) \subset \mathcal{C}^0([a,b]),$$
where $\mathcal{C}^0([a,b])$ is
the space of continuous functions on $[a,b]$
(see \cite{Adams,ABM} for example)  and
$$\forall u \in W^{1,1}(a,b) , \forall y\in \set \qquad \|u\|_{L^\infty} \le |u(y)| + \|u'\|_{L^1}~.$$
\subsection{Gagliardo's fractional Sobolev Spaces}
This section is devoted to recalling the classical definition of fractional Sobolev spaces in the sense of Gagliardo:

\begin{defn}[{\bf Gagliardo's spaces}] Let $s\in(0,1)$. For any  $p\in [1,+\infty)$ we define the following space:

\begin{equation}\label{Gagliardo_def}
\wspG = \left\{u\in \lp\;:\; \frac{|u(x)-u(y)|}{|x-y|^{\frac 1 p + s}}\in L^p(\set \times \set) \right\}\,.
\end{equation}
This is  a Banach space endowed with the norm
$$\|u\|_{\wspG} = \left[ \intset |u(x)|^p d x \;+\; \intset \intset  \frac{|u(x)-u(y)|^p}{|x-y|^{1 + sp}}d x\,d y \right]^{\frac 1 p}\,.$$
\end{defn}
$\wspG$ is an intermediate space between $\lp$ and $\mathrm{W}^{1,p}(a,b )$ and the term
\begin{equation}\label{normG} [u]_{\wspG} = \left[  \intset \intset  \frac{|u(x)-u(y)|^p}{|x-y|^{1 + sp}}d x\,d y \right]^{\frac 1 p}\end{equation}
is the so-called Gagliardo semi-norm of u. We have in particular \cite{palatucci12} that
$$\wspG \subseteq \mathrm{W}_G^{s',p}(a,b ) \quad \quad \forall\quad  0<s'<s\le 1\,.$$
If $s= m+\tau >1$ with $m\in \mathbb{N}, \,\tau\in [0,1[$, such a definition can be generalized to higher orders by setting
$$\wspG = \{u\in \mathrm{W}^{m,p}(a,b ) \;:\;D^m u \in \mathrm{W}^{\tau,p}(a,b) \}.$$
\\
This point of view is quite related to interpolation theory   (see
 \cite{Adams,Aronszajn, demengel12, Gagliardo, Slobodeckij,Tartar} for example). Notice that we are not only interested in the Hilbert framework, so that  a Fourier definition could be useful only if $p=2$.
 There is a huge literature concerning these fractional differentiation methods that we cannot mention here.
\subsection{Fractional integration and differentiation theory}
Another point of view to deal with fractional derivatives is the  one  we describe in the sequel: the  generic definition is the Riemann-Liouville one, though there are many variants that we do not consider. The point of view is different from the  Gagliardo one. Precisely, we get a pointwise definition  of derivatives by using fractional integrals while the Gagliardo's fractional Sobolev Spaces are defined by interpolation and  global approach. As we already mentioned it, we decided to focus on this second type which seems more suitable with respect to applications.
\subsubsection{Fractional integrals}

In what follows $[a,b]$ is a non empty interval of $\R$. We start by defining the fractional integral for $L^1$-functions:

\begin{defn}\label{frac-int}
Let $u\in L^1([a,b])$. For every $s\in (0,1)$ we define the left-side and right-side Riemann-Liouville fractional integrals, by setting respectively
$$I_{a+}^s[u] : x \ge a \mapsto  \frac{1}{\Gamma(s)}  \int_a^x   \frac{u(t)}{(x-t)^{1-s}} d t  ,$$
$$I_{b-}^s[u] : x \le  b \mapsto \frac{1}{\Gamma(s)}  \int_x^b  \frac{u(t)}{(t-x)^{1-s}} d t .$$
\end{defn}
Here $\Gamma$ stands for the classical  Gamma function \cite{abram}.
The properties of left-side and right-side integrals are similar and we prove in the following the main results in the case of the left-side integral $I_{a+}^s[u]$.
The fractional integration theory has been extensively  studied in \cite{samko}. Next proposition recall the main properties of the fractional integral (see Theorem 2.5 and 2.6 in \cite{samko}):


\begin{prop} For any $s\in (0,1)$, the following properties hold true.
\begin{itemize}
\item[(i)]The fractional integral is a continuous operator on the Lebesgue spaces:
\begin{equation}\label{frac-int-cont}
\begin{array}{c}
I_{a+}^s : \lp \rightarrow \lp\,\;\; p\geq 1 \vspace{0.2cm}\\

\|I_{a+}^s u\|_{\lp}\leq C(a,b,s) \| u\|_{\lp}\,.
\end{array}
\end{equation}

\item[(ii)] For every $u\in \lp$, with $p\geq 1$, we have
\begin{equation}\label{frac-int-conv-s}
\underset{s\rightarrow 0_+}{\lim} \|I_{a+}^s u - u\|_{\lp} = 0\,.
\end{equation}
\end{itemize}
\end{prop}
Next theorem concern the mapping properties of fractional integral on Lebesgue and H\"older spaces:
\begin{prop}\label{frac-int-lp} For any $s\in (0,1)$,  we get
\begin{enumerate}
\item[(i)]  $I_{a+}^s$ is a continuous operator from $\lp$ into $L^r(a,b)$ for every $p\in [1,1/s)$ and $r\in [1,p/(1-sp))$;
\item[(ii)] For every $p>1/s$ the fractional integral $I_{a+}^s$ is a continuous operator from $\lp$ into $\mathcal{C}^{0,s-\frac 1 p}(a,b)$;

\item[(iii)] For $p=1/s$  the fractional integral $I_{a+}^s$ is a continuous operator from $\lp$ into $L^r(a,b)$ with $r\in [1,\infty)$;

\item[(iv)]the fractional integral $I_{a+}^s$ is a continuous operator from $\linf$ into $\mathcal{C}^{0,s}(a,b)$.
\end{enumerate}
\end{prop}
Here $\mathcal{C}^{0,s}(a,b)$ denotes the space of H\"older (continuous) functions  of order $s$.
For the proofs we refer to \cite{samko} ($(i)$: proof of Theorem 3.5 p.66, $(ii)$:Theorem 3.6 p.67, $(iii)$: paragraph 3.3 p.91, $(iv)$: Corollary 2 p.56). The  previous proposition shows  that the fractional integration improves the  function regularity.

\begin{rem}[{\bf Fractional integral of BV-functions}]
We point out that to ensure the   H\"older-regularity of the  fractional integral we need to work with  $L^p$-functions with $p>1$. The case $p=1$ is not covered from the previous proposition. \\
However, the point $(iv)$ guarantees such a regularity for bounded functions,  which helps to study an important subset of $\lone$, namely $\BV$ (see section \ref{comparBV}). Indeed,  in dimension  one, every function of bounded variation is bounded, so we get
\begin{equation}\label{frac-int-bv}
I_{a+}^s(\BV)\subset \mathcal{C}^{0,s}(a,b) \quad \quad \forall\; s\in (0,1)\,.
\end{equation}
\end{rem}
Next result states a stronger result for H\"older continuous functions. For the proof we  refer to \cite{samko} again (Theorem 3.4 p.65, Lemma 13.2 p. 240, and Theorem 13.13 p. 238).

\begin{thm}\label{frac-int-hold}
Let $s,\,\alpha\in(0,1)$ 
be such that $s+\alpha < 1$.\\ Then the fractional integral $I_{a+}^s$ is an isomorphism between $\mathcal{C}^{0,\alpha}(a,b)$ and $\mathcal{C}^{0,s+\alpha}(a,b)$: 
\\ there exist two positive constants $C,\,D$, such that,
for every $f\in  \mathcal{C}^{0,s+\alpha}(a,b)$, there exists unique $u\in \mathcal{C}^{0,s}(a,b) $ such that $f=I_{a+}^s[u]$ and
$$C\|f\|_{ \mathcal{C}^{0,s+\alpha}(a,b)} \leq \|u\|_{  \mathcal{C}^{0,s}(a,b)} \leq D\|f\|_{ \mathcal{C}^{0,s+\alpha}(a,b)}\,. $$
\end{thm}

\subsubsection{ Fractional derivatives and representability}
There are several different definitions of fractional derivatives. We recall  next the definition of Riemann-Liouville and Marchaud derivatives and  refer to \cite{samko, samko_exmp} for a deeper analysis of the fractional differentiation theory.

\begin{defn}[{\bf Riemann-Liouville fractional derivative}]\label{defRLderiv}
Let $u\in L^1(a,b)$ and \linebreak $n-1\leq s < n$ ($n$ integer). The left Riemann-Liouville derivative of $u$ at $x\in \set$ is defined by

\begin{equation}\label{RL}
 D_{a+}^s u(x) = \frac{d^n}{d x^n} I_{a+}^{n-s} [u] (x) =   \frac{1}{\Gamma(n-s)} \frac{d^n}{d x^n} \int_a^x   \frac{u(t)}{(x-t)^{s-n+1}} d t\,
 \end{equation}
when it exists.\\
If such a derivative  exists  at $x$ for $s=0$, $n=1$, then it coincides with the function $u$  at $x$.
\\
 Similarly, we may define the right Riemann-Liouville derivative of $u$ at $x\in \set$ as
 \begin{equation}\label{RLr}
 D_{b-}^s u(x) = \frac{d^n}{d x^n} I_{b-}^{n-s} [u] (x) = \frac{1}{\Gamma(n-s)} \frac{d^n}{d x^n} \int_x^b   \frac{u(t)}{(t-x)^{s-n+1}} d t\,
 \end{equation}
 if the last term exists
\end{defn}

\begin{rem} We point out that  $D_{b-}^s$ is the adjoint operator of  $D_{a+}^s$ in the sense to be  precised by Theorem \ref{int-parts}.
\end{rem}
\begin{exa}[{\bf Power function}]\label{power} We consider the function $u(x)=x^k$ ($k\geq 0$) on $[0,1]$ and $a=0$. Then for every $s\in [n-1, n)$ and any $x \in (0,1) $ the fractional derivative at $x$ is defined as
$$
\begin{array}{ll}
 D_{0+}^s  x^k & = \displaystyle{\frac{1}{\Gamma(n-s)} \frac{d^n}{d x^n} \int_0^x   t^k(x-t)^{n-s-1} d t\, }\\[0.3cm]&= \displaystyle{ \frac{1}{\Gamma(n-s)} \frac{d^n}{d x^n}\left[ x^{n+k-s} \int_0^1 (1-v)^{n-s-1}v ^k  d v\,\right] }  \\[0.4cm]
& =
 \displaystyle{ \frac{1}{\Gamma(n-s)} \frac{d^n}{d x^n}\left[ x^{n+k-s}
 B(n-s,k+1)  \,\right] }
 \\[0.4cm]
& =\displaystyle{\frac{1}{\Gamma(n-s)} \frac{\Gamma(k+1)\Gamma(n-s)}{\Gamma(k+1+n-s)} \frac{d^n}{d x^n} x^{n+k-s}} =  \displaystyle{\frac{\Gamma(k+1)}{\Gamma(k-s+1)}x^{k-s}}
\end{array}
$$
where we exploited the Beta Euler function
$B(\nu,\mu)=\int_0^1t^{\nu-1}(1-t)^{\mu-1}=\frac {\Gamma(\nu)\Gamma(\mu)}{\Gamma(\nu+\mu)}$
and
$$\frac{d^n}{d x^n} x^{\tau}=\frac{\Gamma(\tau+1)}{\Gamma(\tau-n+1)}x^{\tau-n}\quad \forall \, \tau\geq 0\,.$$
If $s$ is a positive  integer  number,  then the fractional derivative $D_{0+}^sx^{k}$ coincides with the classical one.\\  If $k=0$, then the fractional derivative $D_{0+}^s$ of the constant function is null if and only if
 $s$ is a strictly positive integer number.
\\
Remark also that
$$D_{0+}^s x^{s-k} = 0 \quad \quad \forall\, s>0 \;, k=1,...,1+[s]\,$$
where $[s]$ denotes the integer part of the real number $s$.
\end{exa}
 Now we focus on the case $n=1$.

\begin{defn}[{\bf Representability}] A function $f\in \lone$ is said to be represented by a fractional integral if $f\in  I_{a+}^s(\lone)$ for some $s\in (0,1)$.
\end{defn}
Theorem 2.1 in \cite{samko} yields the following important result:

\begin{thm}\label{thm-l1}{\bf [$L^1$-representability]} Assume $f\in \lone$.\\ Then $f\in  I_{a+}^s(\lone)$ for some $s\in (0,1)$ if and only if
$$I_{a+}^{1-s} [f] \in W^{1,1}(a,b)\quad \mbox{and} \quad I_{a+}^{1-s} [f] (a) =0\,.$$
Moreover, in the affirmative  case there exists $u\in \lone$ such that $f= I_{a+}^{s} [u]$, and it holds
\begin{equation}\label{inv-l1}
 u= D_{a+}^s   f  \,.
\end{equation}
\end{thm}
As an immediate consequence of Theorem \ref{thm-l1}, we have the following result:
\begin{coro}\label{coro-l1} Assume  $ s\in (0,1)$;   \begin{equation}\label{DI}
\forall u\in \lone \qquad D_{a+}^s  I_{a+}^s [u] = u ,
\end{equation}
and
\begin{equation}\label{ID}
\forall u\in I_{a+}^s(\lone) \qquad I_{a+}^sD_{a+}^s  u  = u .
\end{equation}
\end{coro}
In the previous Corollary, \eqref{DI} proves that fractional differentiation  can be seen as the inverse operator of the fractional   integration.
\\
The converse is not true in general:  a counterexample is given by the power function $x^{s-k}$ ($k=1,...,1+[s]$) whose $s$-fractional derivative  is null if and only if $k=0$ and $s>0$  (Example \ref{power}).  This is similar  to the  classical integration and differentiation theories where the integral of $u'$ differs from $u$ for a constant.
However, according to \eqref{ID}, for every function that  can be represented as a fractional integral, the fractional    integration acts as the reciprocal operator of the fractional differentiation.

Using Theorem \ref{thm-l1}, it is easy to verify that the power function $x^{s-1}$  is not represented by a fractional integral. In fact $I_{a+}^{1-s} [x^{s-k}] \equiv \Gamma(s-1)$ so that it belongs to $W^{1,1}(a,b)$ but it does not verify $I_{a+}^{1-s} [x^{s-k}] (0) =0$.

Next proposition follows from Theorems \ref{frac-int-hold} and \ref{thm-l1}:

\begin{prop}\label{exist-der-hold}
Let $\alpha \in (0,1)$ and $u\in \mathcal{C}^{0,\alpha}(a,b)$. Then $D_{a+}^s u$ exists for every $ s \in [0,\alpha)$ and $D_{a+}^s u \in \mathcal{C}^{0,\alpha-s}(a,b)$.
\end{prop}

 Let us give a comment about the relationship  between $D_{a+}^s$ and $I_{a+}^s$.
There are two kind of results  on  the fractional integral that can be very useful
  to study the properties of the fractional derivative :
\begin{itemize}
\item
The first-one is a representability result (for instance Theorem \ref{thm-l1}) that gives conditions for   a function  $f$ to be  represented as the fractional integral of an other function $u$. This is quite  important because it allows to easily prove that $f$ admits a fractional derivative  with \eqref{inv-l1}.

However,  the representability of a function $u$ (i.e., $u=I^s_{0+}\varphi$) is only a sufficient condition to get  the existence  of the derivative. The power function and the Heaviside function give two examples of functions that are not representable and whose fractional derivative exists (see  previous discussion and Example \ref{charact}).

\item The second kind  of result  	are  embedding results, as \eqref{frac-int-bv}, that give some  informations  on  the regularity of the fractional integral to get  the Riemann-Liouville fractional derivative existence.
\end{itemize}

\subsubsection{Marchaud derivative and representability for $p\in (1,\infty)$}

The representability result given by  Theorem \ref{thm-l1} can be improved by characterizing the set  of the functions $u\in L^p(a,b)$ ($p>1$) represented by an other $L^p$-function
$$u=I_{a+}^s[f]\quad f \in L^p(a,b)\;\; s\in (0,1)\,.$$
In order to address this issue  we need to introduce a slightly different definition/notion  of fractional derivative.
\\
According to \cite{samko} p. 110, we note that, for $\mathcal{C}^1$-functions and  every $s\in (0,1)$,  the use of integration by parts  gives
\begin{equation}\label{Marchaud-smooth}
 D_{a+}^s u(x) = \frac{u(x)}{\Gamma(1-s)(x-a)^s} + \frac{s}{\Gamma(1-s)}\int_a^x \frac{u(x)-u(t)}{(x-t)^{1+s}}\,dt \quad \forall \,x\in (a,b]\,.
 \end{equation}
The \textit{Marchaud fractional derivative} is defined as  the right  term  of  \eqref{Marchaud-smooth}. To  extend this setting   to   non-smooth functions we need to define the integral by a limit, which leads to the following definition:

\begin{defn}[{\bf Marchaud fractional derivative}]
Let $u\in \lone$ and $s\in (0,1)$. The left-side  Marchaud derivative of $u$ at $x\in (a,b]$ is defined by
$$\DD_{a+}^s u(x) = \underset{\varepsilon\rightarrow 0}{\lim}\;\DD_{a+, \varepsilon}^s u(x)$$
with
$$\DD_{a+, \varepsilon}^s u(x) = \frac{u(x)}{\Gamma(1-s)(x-a)^s} + \frac{s}{\Gamma(1-s)} \,\psi_\varepsilon(x)$$
and
\begin{equation}\label{psi_e}
\psi_\varepsilon(x)=
\left\{
\begin{array}{ll}
\displaystyle{\int_{a}^{x-\varepsilon} \frac{u(x)-u(t)}{(x-t)^{1+s}}\,dt}& \mbox{if}\; x \geq  a+\varepsilon\,,\\
\displaystyle{\int_{a}^{x-\varepsilon} \frac{u(x)}{(x-t)^{1+s}}\,dt} & \mbox{if} \;\;a\leq x\leq a+\varepsilon\,.
\end{array}\right.
\end{equation}
\end{defn}
Note that the definition of $\psi_\varepsilon$ for $a\leq x\leq a+\varepsilon$ is obtained by continuing the function $u$ by zero beyond the interval $[a,b]$.
The  passage to the limit depends  on the functional space we are working with. In our case we will consider $L^p$-functions and the limit will be defined in the sense of the $L^p$-strong topology.

We remark that  such a derivative is not   defined at $x=a$ and that  a necessary condition for the derivative to exist  is    $u(a)=0$

The right-side derivative can be defined similarly by using the integral between $x$ and $b$. In the following we state the main results about Marchaud differentiation for the left-side derivative, but  similar results can be obtained for the right-side one.

As expected,  if $u$ is a  $\mathcal{C}^1$-function, the Marchaud and Riemann-Liouville derivative coincide for every $s\in (0,1)$, and their expression is given by \eqref{Marchaud-smooth}:
$$ \forall\, u\in \Cone\qquad  D_{a+}^s u(x) = \DD_{a+}^s u(x)\quad \forall \,x\in(a,b].$$

Next results  generalize Theorem \ref{thm-l1} and Corollary \ref{coro-l1} :
\begin{thm}\label{thm-rep-lp} Let  $s\in (0,1)$ and $p\in (1,\infty)$.  For any $f\in \lp$, we get  $f=I_{a+}^s[u]$ with $u\in \lp$ if and only if  the limit of the family $\{\psi_\varepsilon\}$ as $\varepsilon \rightarrow 0$,  where $\psi_\varepsilon$ is defined in \eqref{psi_e}   exists (for the    $L^p$  norm topology).
\end{thm}

\begin{thm}\label{thm-lp}
Let be  $f=I_{a+}^s[ u]$ where  $u\in \lp$ with $p\geq 1$ and $s\in (0,1)$. Then
$$\DD_{a+}^s f=u\,.$$
\end{thm}
Proofs can be found in \cite{samko} (Theorem 13.2 p. 229  and Theorem 13.1 p. 227, respectively).

\begin{rem}[{\bf Marchaud vs Riemann-Liouville derivative}]

We point out that, for every $s\in(0,1)$, we have
$$\forall\, u\in I_{a+}^s(\lone),~ \forall \,x\in(a,b] \quad  D_{a+}^s u(x) = \DD_{a+}^s u(x),$$
because of Theorems \ref{thm-l1} and \ref{thm-lp}.
With Theorem \ref{frac-int-hold}, this implies in particular that
$$ \forall\, u\in \mathcal{C}^{0,s+\alpha}(a,b), ~ \forall  x\in (a,b]\qquad D_{a+}^s u(x) = \DD_{a+}^s u(x),$$
if $s+\alpha<1$. This is a useful result in order to study the fractional derivative because the Marchaud derivative it is easier to handle.
\end{rem}
\section{Riemann-Liouville  Fractional Sobolev  space ($p=1$)}

In this section, we define the   Sobolev spaces associated to the Riemann-Liouville fractional derivative   for $p=1$. The case $p=1$ is of particular interest since we aim to  study the relationship between these spaces and the spaces of functions of bounded variation.

A first possible definition could be given by the following set for  $ s\in(0,1)$:
$$\{u\in \lone\;|\; D_{a+}^su \in \lone\}~, $$
which contains all the  $L^1$-functions  such that the   Riemann-Liouville fractional derivative  or order $s$ for a given $s\in (0,1)$ belongs to $L^1$.

We noticed that if the Riemann-Liouville fractional derivative of $u$  exists  (for some $s$), then   $I_{a+}^s[u]$ is differentiable almost everywhere. However, we have no  information  on  the  differential properties of the fractional integral.  These   differential properties are not completely described by the  pointwise derivative though it exists a.e. (it could be for example an $SBV$-function).
This shows that the previous definition is not suitable to  obtain a generalized integration by parts formula.


Therefore,  to develop a satisfactory theory of fractional Sobolev spaces  we use a more suitable definition  in the next section.

\subsection{Riemann-Liouville  Fractional Sobolev  spaces}

Following \cite{Bourdin15,Idczak13}  where these spaces are denoted  $AC^{s,1}_{a^+}$,  we may define the Riemann-Liouville  Fractional Sobolev spaces  as follows:

\begin{defn} Let $s\in[0,1)$. We denote by
$$ W^{s,1}_{RL,a+} := \{ u \in \lone \;|\; I_{a+}^{1-s}[u]\in W^{1,1}(a,b)\;  \}.$$
A similar space can be defined by using the right-side fractional integral.
\end{defn}
Note that this definition does not mean that $u$ is representable but  its fractional integral  $f= I_{a+}^{1-s}u$ is representable.  We shall describe the representable functions of  $W^{s,1}_{RL,a+}$ next (it is related to their trace).  \\
We can make the connection with what would be the \textit{natural} definition
\begin{prop}  Let $s\in[0,1)$.  Then
$$ W^{s,1}_{RL,a+} \subset  \{ u \in \lone~|~  D_{a+}^su \in \lone~\}.$$
\end{prop}
\proof Let be  $u\in  W^{s,1}_{RL,a+}$ $\subset  \lone $. With definition \ref{defRLderiv} (with $n=1$)  we get
 $\displaystyle{D_{a+}^s u = \frac{dI_{a+}^{1-s} [u] }{d x}  }$.
  In addition, Theorem 1 p. 235 in \cite{EG}, yields that
 if $I_{a+}^{1-s}[u]\in W^{1,1}(a,b)$ its Fr\'echet  derivative exists a.e. and coincides with its weak derivative. \\This means that if $u\in W^{s,1}_{RL,a+}(a,b) $ then the fractional derivative $D_{a+}^s u$ exists a.e. and belongs to $\lone$.
\proofend
\begin{rem}  If  $I_{a+}^s[u]$ is differentiable at every point and its derivative is a $L^1$-function, then $u\in W^{s,1}_{RL,a+} (a,b)$. This is due to the fact that if a  function  is differentiable at every point and its derivative is $L^1$, then it belongs to $W^{1,1}(a,b)$.
\end{rem}
In the sequel, when  there is no ambiguity, we omit the index ``RL,a+''  in the notation.
\\
Analogously to the usual theory we introduce the following space:

\begin{defn} Let $s\in(0,1)$ and $p\in [1,\infty)$. We denote by
$$\wszero := \{ u \in \lone \;|\; I_{a+}^{1-s}[u]\in W^{1,1}_{0,a}(a,b)\}.$$
where  $W^{1,1}_{0,a}(a,b)=\{v\in W^{1,1}(a,b)\;|\; v(a)=0\}$.
\end{defn}
Note that $u \in \wszero$    does not mean that  $u(a) = 0$ (it may  even not be defined). \\
The following result  is a direct consequence  of Theorem \ref{thm-l1}:

\begin{thm}\label{rep_W0} Let $s\in[0,1)$. Then $u\in \wszero$ if and only if $u$ is $L^1$-representable.
\end{thm}
%

\subsection{Main properties}

Before performing comparisons between these fractional Sobolev spaces and the   spaces of bounded variation functions, we investigate their basic properties.

\begin{thm} For any $s \in (0,1)$,
the Riemann-Liouville fractional Sobolev space  $\wsoneRL$ is a Banach space endowed with the norm
$$ \|u\|_{\wsoneRL}:= \|u\|_{\lone} + \|I_a^{1-s}[ u ]\|_{W^{1,1}(a,b)}.$$
\end{thm}

\proof It is easy to verify that $\|\cdot\|_{\wsoneRL}$ is a norm so that we  have just to prove the completeness. Let $(u_n)_{n \in \N}$ be a Cauchy sequence in $\wsoneRL$ which implies that $(u_n)_{n \in \N}$  and
$(I_{a+}^{1-s} [u_n])_{n \in \N}$ are Cauchy sequences in $L^1(a,b)$ and $W^{1,1}(a,b)$, respectively. Then there exists $u\in L^1(a,b)$ and $v\in  W^{1,1}(a,b)$ such that
$$u_n\,\overset{L^1}{\rightarrow} \, u\;, \quad \quad I_{a+}^{1-s} [u_n]\,\overset{W^{1,1}}{\rightarrow} \, v$$
Because of the definition we have $I_{a+}^{1-s}[u_n]\in W^{1,1}(a,b)$, and,
 because of \eqref{frac-int-cont}, as $u_n\rightarrow u$ strongly in $L^1$, we  have $I_{a+}^{1-s}[u_n]\rightarrow I_{a+}^{1-s}[u]$
strongly in $L^1$ as well. This proves that $I_{a+}^{1-s}[u]=v$.
\\
Now, we have to prove the (strong) $L^1$ convergence of the first distributional  derivative of  $(I_{a+}^{1-s}[u_n])'$ toward  the   derivative of  $I_{a+}^{1-s}[u]$ that is : $ v'= (I_{a+}^{1-s}[u])'$, where $v'$ denotes the weak derivative of $v$.
For every  $\varphi\in C^\infty_c(a,b)$ we get
$$\int_a^b \varphi\, (I_{a+}^{1-s}[u_n])' =-  \int_a^b   \varphi'  I_{a+}^{1-s}[u_n]  $$
and, by taking the limit, we obtain
$$\int_a^b \varphi\, v' = -  \int_a^b  \varphi'  I_{a+}^{1-s}[u].$$
This proves that $I_{a+}^{1-s}[u]$ is a $W^{1,1}$-function and its weak derivative is $v'$. This implies that $I_{a+}^{1-s} [u_n]\,\overset{W^{1,1}}{\rightarrow}I_{a+}^{1-s} [u]$ and  this gives the result.
  \proofend
 
 An immediate consequence of this result is the following theorem
 \begin{thm} $\wszero$ equipped with the norm
 $$ \|u\|_{\wszero}:= \|u\|_{\lone} + \| D_{a+}^{s} u \|_{L^{1}(a,b)}$$
  is a Banach space.
\end{thm}
\proof The space $\wszero$ is complete with respect to the norm $ \|\cdot\|_{\wsoneRL}$. This is done  as in the previous proof  since  $W^{1,1}_{0,a}(a,b)=\{u\in W^{1,1}(a,b)\;|\; u(a)=0\}$ is a Banach space.
\\
Moreover,  with Poincar\'e's inequality  for spaces $W^{1,1}_{0,a}(a,b)$, the norms $ \|\cdot\|_{\wsoneRL}$ and $\|\cdot\|_{\wszero}$ are equivalent.
\proofend
  \begin{thm}
Assume $0<s<s'<1$  and consider  $u\in  I^{s'}_{a+}(L^1(a,b))$.  Then  $ u \in W_0^{s,1}(a,b)$ and
  $$\|u\|_{W^{s,1}(a,b)}\leq C_{s,s'}\|u\|_{W^{s',1}(a,b)}\,.$$
   \end{thm}
   \proof  Let be $u\in I^{s'}_{a+}(L^1(a,b))$.  Then $u=I^{s'}_{a+}[f] $ with $f\in L^1(a,b)$. So  with  Theorem 2.5 p.46 in \cite{samko} we get  $u=I^{s'}_{a+}[f]=I^{s}_{a+}[I^{s'-s}_{a+}[f]]$. Therefore  $u\in I^{s}_{a+}(L^1(a,b))$ for every $s\in (0,s']$.    \\
  As $u$ is represented by a fractional integral of a $L^1$-function, then  with Theorem  \ref{rep_W0},  $ u \in W_0^{s,1}(a,b)$ for  every $s\in (0,s']$.  As $f \in L^1(a,b)$, Corollary  \ref{coro-l1}  gives
  $$ \|D^{s'}_{a+} u \|_{L^1} =  \|D^{s'}_{a+} I_a^{s'} [f] \|_{L^1}  = \| f \|_{L^1}~.$$
  Similarly, $u$ is represented by  $g=  I^{s'-s}_{a+}[f] $ and
   $$ \|D^{s}_{a+} u \|_{L^1} = \| g \|_{L^1} =  \| I^{s'-s}_{a+}[f]  \| _{L^1}~.$$
   Moreover,  with the continuity of the fractional integral operator we get
  $$ \| I^{s'-s}_{a+}[f]  \| _{L^1} \le  C_{s,s'} \| f \|_{L^1 },$$
  where $ C_{s,s'}$ is a generic constant depending on $s$ et $s'$.
  Finally
  $$ \|D^{s}_{a+} u \|_{L^1}  \le  C_{s,s'}\|D^{s'}_{a+} u \|_{L^1} ~.$$
  This proves the result.
%
%
    \proofend
Next theorem  gives a relationship between Riemann-Liouville Sobolev spaces and Gagliardo Sobolev spaces:

\begin{thm} Let  be $s ,s' \in (0,1)$ such that $s'>s$.  Then
$$ W^{s',1}_G(a,b) \cap I^s_{a+}(\lone) \subset  W^{s,1}_{RL,a+}(a,b)$$
with continuous injection. More precisely,
   $$\forall u\in W^{s',1}_G(a,b)\cap I^s_{a+}(\lone)  \qquad \|u\|_{W^{s,1}_{RL,a+}(a,b)}\leq C \|u\|_{W^{s',1}_{G}(a,b)}\,.$$
  \end{thm}
  \proof Let us choose $s \in (0,1) $ and $s'>s$ (in $(0,1)$).
Let be  $u\in I^s_{a+}(\lone)$. It is represented by a fractional integral of a $L^1$-function, so  its Riemann-Liouville and Marchaud derivative coincide.  The norm of the fractional derivative can be evaluated by the Marchaud derivative.  Recall that
if $a\le x\le a+\varepsilon$ we get
$$\DD_{a+, \varepsilon}^s u(x)= \frac{u(x)}{ \varepsilon^s\Gamma(1-s)},$$
and  if $  x\ge a+\varepsilon$
$$\DD_{a+, \varepsilon}^s u(x)=  \frac{1}{\Gamma(1-s)} \frac{ u(x)}{(x-a)^s}+   \frac{s}{\Gamma(1-s)} \int_{a}^{x-\varepsilon} \frac{u(x)}{(x-t)^{1+s}}\,dt.$$
   \begin{align*}
 \|\DD_{a+, \varepsilon}^s u\|_{L^1(a+\varepsilon,b)}\leq & \frac{1}{\Gamma(1-s)}\int_{a+\varepsilon}^b \frac{|u(x)|}{(x-a)^s} ds
 +   \frac{s}{\Gamma(1-s)}\int_{a+\varepsilon}^b\int_{a}^{x-\varepsilon} \frac{|u(x)-u(t)|}{|x-t|^{1+s}} \, dt \,dx\\
 \leq &   \frac{1}{\Gamma(1-s)}\int_{a+\varepsilon}^b \frac{|u(x)|}{(x-a)^s} ds   +  \frac{s}{\Gamma(1-s) }[u]_{W^{s,1}_{G}(a,b)} .
 \end{align*}
where  the Gagliardo semi-norm  $[u]_{W^{s,1}_{G}(a,b)}$  is  given by \eqref{normG}. Moreover
 $$ \|\DD_{a+, \varepsilon}^s u\|_{L^1(a,a+\varepsilon)}\leq  \frac{1}{ \varepsilon^s\Gamma(1-s)}\int_a^{a+\varepsilon} |u(x)| ds; $$
  finally

\begin{equation}\label{WG2}
 \|\DD_{a+, \varepsilon}^s u\|_{L^1(a,b)}\leq  \frac{1}{\Gamma(1-s)}\left (\int_{a+\varepsilon}^b \frac{|u(x)|}{(x-a)^s} ds  + \varepsilon^{-s}\int_a^{a+\varepsilon} |u(x)| ds +  s [u]_{W^{s,1}_{G}(a,b)}\right).
 \end{equation}
 Now, we know (\cite{palatucci12} - section 6.   for example) that  $$W^{s',1}_{G}(a,b) \subset L^p(a,b)$$
 with continuous injection and  $p \in \displaystyle{[1, \frac{1}{1-s'}})$. As     $u \in  W^{s',1}_G(a,b) $, then  $ u \in L^p(a,b)$ where $p$ can be chosen such as $  \frac{1}{1-s}< p <  \frac{1}{1-s'}$   since that $s' > s$. Let us call $p^*$ the conjugate exponent and apply H\"older inequalities to relation \ref{WG2}. Note that $p^*$ satisfies $1-sp^* >0$.   We get

 \begin{align*}
  \|\DD_{a+, \varepsilon}^s u\|_{L^1(a,b)}\leq & ~\frac{ \|u\|_{L^p(a+\varepsilon,b)} }{\Gamma(1-s)}
 \left [ \frac{(b-a)^{1-sp^*} - \varepsilon^{1-sp^*}}{ 1-sp^*}\right]^{1/p^*}   \\
 &
  + \frac{ \|u\|_{L^p(a, a+\varepsilon)} }{\Gamma(1-s)}   \varepsilon^{(1-sp^*)/p^*}  \\
& +  \frac{s}{\Gamma(1-s)}   [u]_{W^{s,1}_{G}(a,b)}
  \\
  \leq & \frac{\|u\|_{L^p}}{ \Gamma(1-s)}
  \left ( \left [ \frac{ (b-a)^{1-sp^*} - \varepsilon^{1-sp^*} }{1-sp^*} \right]^{1/p^*}
  + \varepsilon^{(1-sp^*)/p^*}\right) \\
 & +  \frac{s}{\Gamma(1-s)}  [u]_{W^{s,1}_{G}(a,b)}~.
 \end{align*}

 Passing to the limit as $\varepsilon \to 0$ gives
 $$
  \|\DD_{a+ }^s u\|_{L^1(a,b)}\leq \frac{\|u\|_{L^p}}{\Gamma(1-s)}  \left( \frac{ (b-a)^{(1-sp^*)} }{   1-sp^* }\right)^{1/p^*}
  +  \frac{s}{\Gamma(1-s)}  [u]_{W^{s,1}_{G}(a,b)}~.
 $$
 As   $ \|u\|_{L^p} \le C  \|u \|_{W^{s',1}_{G}(a,b)}$  and
 $$ [u]_{W^{s,1}_{G}(a,b)} \le  [u]_{W^{s',1}_{G}(a,b)}\le   \|u \|_{W^{s',1}_{G}(a,b)}$$ we finally get
 $$  \|\DD_{a+ }^s u\|_{L^1(a,b)}\leq  C(s,s',a,b) \|u\|_{W^{s',1}_{G}(a,b)}~.
 $$
 This ends the proof.
  \proofend

  Note that $s=s'$ is the critical case in the above theorem. We cannot handle this case with the same arguments.

 \subsection{Integration by parts and relationship with $W^{1,1}(a,b)$ }

We have the following integration by parts formula (\cite{Bourdin15} Theorem 2).

 \begin{thm}If $\displaystyle{0 \le \frac{1}{p} < s < 1 }$ and   $\displaystyle{0 \le \frac{1}{r} < s < 1 }$, then
 for every $u \in \wsp$ and $v \in \wsm$ we get
 \begin{align} \label{IPP}
 \int_a^b   (D_{a+}^s u ) (t) v(t) \, dt =  \int_a^b   (D_{b-}^s v ) (t) u(t) \, dt +
 u(b) \,I_{b-}^{1-s} [v](b) -  I_{a+}^{1-s} [u](a) \,v(a)~.
 \end{align}
 \end{thm}
 Here
 $$ W^{s,1}_{RL,b-} := \{ u \in \lone \;|\; I_{b-}^{1-s}[u]\in W^{1,1}(a,b)\;  \}.$$
 This is a generalization of the following proposition proving that  $D_{a+}^s$ and $D^s_{b-}$ are adjoints operators on the set of functions represented by fractional integral (see Corollary 2 p. 46 in \cite{samko})
\begin{thm}\label{int-parts} Let $u\in I^s_{a+}(L^p(a,b))$ and  $v\in I^s_{b-}(L^q(a,b))$  with $p,q\in [1,+\infty[$, $1/p+1/q\leq 1+s$ and $p,q\neq 1$ if $1/p+1/q= 1+s$. Then
 $$
  \int_a^b   (D_{a+}^s u ) (t) v(t) \, dt =  \int_a^b   (D_{b-}^s v ) (t) u(t) \, dt\,.
$$
\end{thm}

\begin{exa}[{\bf Smooth functions}]\label{smooth} Set $[a,b]= [0,1]$. For every $u\in C^\infty([0,1],\R)$ and for every $s\in[0,1)$ we get for any $x \neq 0$:

\begin{equation*}
\label{smooth2}
\begin{array}{ll}
 D_{0+}^s u(x)  & = \displaystyle{\frac{1}{\Gamma(1-s)} \frac{d}{d x} \int_0^x   u(t)(x-t)^{-s} d t }\vspace{0.2cm}\\
 &=\displaystyle{  \frac{1}{(1-s)\Gamma(1-s)}   \frac{d}{d x}\left[u(0)x^{1-s}+ \int_0^x   u'(t)(x-t)^{1-s} d t\right] }  \vspace{0.2cm}\\
 & =\displaystyle{ \frac{1}{\Gamma(1-s)}\left[ u(0)x^{-s} + \int_0^x   u'(t)(x-t)^{-s} d t\right] }\vspace{0.2cm}\\
 &=\, \displaystyle{\frac{1}{\Gamma(1-s)} u(0)x^{-s}+ \frac{1}{\Gamma(2-s)} u'(0)x^{1-s}+\frac{1}{\Gamma(2-s)}\int_0^x   u''(t)(x-t)^{1-s} d t}\, \,
 \end{array}
 \end{equation*}
since
$\Gamma(z+1) = z\Gamma(z)$ for every  $z >0$.
Note that, in order to make such a derivative well defined at $x=0$, we have to suppose $u(0)=0$.
\\
As $ x \mapsto x^{-s}$ belongs to $L^1(0,1)$ and $u''$ is bounded and $s\in[0,1)$, we get that $D_{0+}^s u$ belongs to $L^1(0,1)$.
Moreover, we have
\begin{align*}
\|u'-D_{0+}^su\|_{L^1([0,1],\R)}= &
\int_0^1\left| u'(x)
 -\frac{u(0)x^{-s}}{\Gamma(1-s)}
- \frac{u'(0)x^{1-s}}{\Gamma(2-s)} -\int_0^x \frac{ u''(t)(x-t)^{1-s}}{\Gamma(2-s)} dt\right|dx\\
\le &\frac{|u(0)|} {(1-s)\Gamma(1-s)}  +
\int_0^1\left| u'(x)
- \frac{u'(0)x^{1-s}}{\Gamma(2-s)} -\int_0^x \frac{ u''(t)(x-t)^{1-s}}{\Gamma(2-s)} dt\right|dx \\
\le &\frac{|u(0)|} {\Gamma(2-s)}  +
\int_0^1\left| u'(x)
- \frac{u'(0)x^{1-s}}{\Gamma(2-s)} -\int_0^x \frac{ u''(t)(x-t)^{1-s}}{\Gamma(2-s)} dt\right|dx
\end{align*}
Since $\Gamma(1)=1$   by the Lebesgue convergence theorem, we get
$$
\|u'-D_{0+}^su\|_{L^1([0,1],\R)}\overset{s\rightarrow 1}
{\longrightarrow}\;
|u(0)|+\int_0^1\left| u'(x) - u'(0)- \int_0^x  u''(t)
d t \right| dx =|u(0)|\,.
$$
Moreover, if $u(0)=0$
$$ \forall x \in [0,1]\qquad  \lim_{s \to 1} D_{0+}^s(u)(x) = u'(x),$$
and not only almost everywhere.
\end{exa}
Now, we compare the usual Sobolev space $ W^{1,1} (a,b)$ and the Riemann-Liouville fractional Sobolev spaces.
\begin{thm} \label{w11} The space   $ W^{1,1}(a,b)$ is continuously embedded in  $\wsp$ for any $s\in [0,1)$ and  $$D_a^s u \;\; {\rightharpoonup}\;\; u'  \mathcal{L}^1+ u(a) \delta_a\,\quad \mbox{as}\,\quad  s\rightarrow 1\,,$$
weakly    in $\mathcal{M}(a,b)$  (the space of Radon measures on $[a,b]$). Here $\delta_a$ is the Dirac measure at $a$ and $\mathcal{L}^1$ denotes the 1D Lebesgue measure.  Moreover,
  if $u(a)=0$ we have
 $$D_{a+}^s u \overset{L^1(a,b)}{\longrightarrow} u'\,\quad \mbox{as}\,\quad  s\rightarrow 1 .$$
\end{thm}
 \noindent\textit{Proof.}  According to Lemma 2.1 p. 32 in \cite{samko} and its corollary, every $u\in W^{1,1}(a,b) $ has a fractional derivative for every $s\in (0,1)$. Fractional integrals and derivatives are given by
\begin{equation}\label{repr-w11} I^{1-s}_{a+} [u](x)= \frac{1}{\Gamma(2-s)} \left[u(a)\,(x-a)^{1-s}+ \int_a^x   u'(t)(x-t)^{1-s} d t\right]\,
\end{equation}
and
\begin{equation}\label{repr-w11-bis}
\begin{array}{lll}
 \forall\,x\in[a,b]&  D^s_{a+}  u(x)&\displaystyle{=\frac{1}{\Gamma(1-s)} \left[u(a)\,(x-a)^{-s}+ \int_a^x   u'(t)(x-t)^{-s} d t\right]}\\
 && \displaystyle{= \frac{u(a)(x-a)^{-s}}{\Gamma(1-s)}+  I^{1-s}_a [u'](x)\,.}\\
 \end{array}
 \end{equation}
As $s\in (0,1)$  with \eqref{frac-int-cont}  we get that the fractional derivative is defined at every point of
$(a,b]$, belongs to $\lone$, and
\begin{align*} \|D^s_{a+}  u\|_{\lone} \leq &\frac{|u(a)| (b-a)^{1-s}}{\Gamma(2-s)} +  \|I^{1-s}_{a+}  [u']\|_{L^1}(a,b) \\
\leq & C  \left( |u(a)|  +  \|  u'\|_{L^1}(a,b) \right ).
\end{align*}
As $W^{1,1}(a,b) $ is continuously embedded in $ \mathcal{C}([a,b])$(see  \cite{Adams,brezis}), we get $|u(a)| \le  \|u\|_\infty \le C \|u\|_{W^{1,1}}$ and
$$\|D^s_{a+}  u\|_{\lone} \leq  C  \|u\|_{W^{1,1}}$$
which proves the continuous embedding. (Here $C$ is a generic constant that only depends on $a,b$ and $s$).
 \\
Now, as $u' \in L^1(a,b) $ then with  relation \eqref{frac-int-conv-s} we get the strong convergence of  $I^{1-s}_a u'$ to $u'$ in $L^1(a,b) $ as $s \to 1$.  Therefore we get the  result in the case where $u(a)=0$ using  \eqref{repr-w11-bis}.
%

 Moreover,  for every $\varphi\in \mathcal{C}^1(a,b)$
\begin{eqnarray*}
\frac{u(a)}{\Gamma(1-s)} \!\! \int_a^b \!\!\!(x-a)^{-s}\varphi(x) dx &\!\!\!\!\!\!=& \!\!\!\!\frac{u(a)}{\Gamma(2-s)} \left[ \varphi(b) (b-a)^{1-s}   - \!\!\!  \int_a^b\!\! (x-a)^{1-s}\varphi'(x) dx \right ]\\ &\!\!\!\!\!\!\overset{s \to 1}{\to}&\!\!\!\! u(a)  \varphi(a)  .\end{eqnarray*}
 We conclude with the density of the  $\mathcal{C}^1$ functions  in the space of continuous functions
 to get  convergence of $D^s_{a+}  u$ to
$ u'\mathcal{L}^1 +  u(a) \delta_a$  in $\mathcal{M}(a,b)$.
\proofend
\\

%
%

\begin{rem} We give an example of a smooth function $u$ with $u(a)\neq 0$ and such that the fractional derivative  does not converge (as $s\rightarrow 1$) to $u'$ strongly in $L^1$.
\\
It suffices to consider $u(x)=1$ on $\set$. Then
$$u'\equiv 0\;,\quad \quad D^s_{a+}  u(x) = \frac{(x-a)^{-s}}{\Gamma(1-s)}$$
and
$$\|D^s_{a+}  u\|_{\lone}\,= \,\frac{(b-a)^{1-s}}{\Gamma(2-s)}\,\overset{s\rightarrow 1}{\longrightarrow}\, 1\,\neq \,\|u'\|_{\lone}\,. $$

\end{rem}

\section{Comparison with BV and SBV}\label{comparBV}
 Let us   recall the  definition and the main properties  of  the  space  of functions of bounded variation  (see \cite{ambrosio,ABM} for example), defined by
$$\bv = \{ u \in L^1(\Omega) \ |\ TV(u) < + \infty\}, $$
where  $\Omega$ is a bounded, open subset of $\R^{d}$ and
\begin{equation}\label{ju}
TV(u)  :=\sup \left \{ \int_\Omega u(x) \, \mbox{div }\xi (x) \, dx  \ |\ \xi \in \mathcal{C}^1_c(\Omega),\ \|\xi\|_\infty\leq 1 \right\}.
\end{equation}
The space $\bv$, endowed with the norm
$ \|u\|_{\bv} = \|u\|_{L^1} + TV(u),$
is a Banach space.
The derivative in the sense of  distributions of every $u\in \bv$ is a bounded Radon measure, denoted $Du$, and  $TV(u)=\int_\W |Du|$ is the total variation of $ u$.
We next recall standard properties of functions  of bounded variation  in the case $d=1$.

\begin{prop}\label{proprecall}
Let  $(a,b)$ be an open subset of $ \R $.
\begin{enumerate}
\item For every $u\in \bv$, the Radon measure $Du$ can be decomposed into
$ Du = \nabla  u \,dx + D^\perp u$,
where $\nabla  u \, dx $ is the absolutely continuous part of $Du$ with respect of the  Lebesgue measure and $D^\perp u$ is the singular part.
\item The mapping $u \mapsto TV(u)$ is lower semi-continuous   from $\bv$ to $\R^+$ for the $L^1(a,b)$ topology.
\item $\bv \subset L^\sigma(a,b)$ with continuous embedding, for  $ \sigma\in [1, \infty]$.
\item $\bv \subset L^\sigma(a,b)$  with compact embedding, for  $\sigma\in [1,\infty  )$ .
\end{enumerate}
\end{prop}

The singular part $D^\perp u$  of the derivative has a  jump part and a Cantor component.  The $SBV(\Omega)$ space ( see \cite{ambrosio} for example)   is the space of functions in
 $BV(\Omega)$ whose derivative has no singular Cantor component.   The functions of
  $SBV(a,b)$  have two components: one is regular and belongs to $W^{1,1}(a,b)$ and the other one is a countable summation of characteristic functions. More precisely, any increasing function  in  $SBV(a,b)$ can be written as
$$u(x)=u(a)+\int_a^xu'(t)dt + \underset{x_k\in J_u}{\sum} p_k\chi_{[x_k,1]}(x) \quad \quad x\in[a,b]\,$$
where $J_u$ denotes the (at most countable) set of jump points of $u$ and $p_k=u^+(x_k)-u^-(x_k)$ denotes the positive jump of $u$ at $x_k$.  This describes all the functions of $SBV(a,b)$ since any $BV$-function can be written as the difference of two increasing functions.

Next example shows that there exists  a $SBV$-function that belongs  to  $\wsp$ for any $s\in (0,1)$.  This confirms the regularizing behavior of the fractional integral operator and represents a preliminary result in order to prove the relationship between $SBV$ function and fractional Sobolev space.

\begin{exa}[{\bf Heaviside function}]\label{charact} Let $u:[0,1]\rightarrow\mathbb{R}$ ($a=0$), $u(x)= \chi_{[\alpha,1]}(x)$ with $\alpha\in (0,1)$. We consider $s\in [0,1)$.  For every $x\in [\alpha,1]$ we get
\begin{equation}\label{charactder}
I^{1-s}_{0+} [\chi_{[\alpha,1]}](x)=
\left\{
\begin{array}{ll}
0 & \mbox{if}\quad x\in[0,\alpha)\\
\displaystyle{\frac{ (x-\alpha)^{1-s}}{\Gamma(2-s)}}& \mbox{if}\quad x\in [\alpha,1] \\
\end{array}\right.
\end{equation}
which proves that $I^{1-s}_{0+}[\chi_{[\alpha,1]}]\in W^{1,1}([0,1])$ so that $ \chi_{[\alpha,1]}\in W^{s,1}_{RL, 0+} ([0,1])$.

The fractional derivative is given by
\begin{equation}\label{charactder-bis}
D^s_{0+} \chi_{[\alpha,1]}(x)=
\left\{
\begin{array}{ll}
0 & \mbox{if}\quad x\in[0,\alpha]\\
\displaystyle{\frac{ (x-\alpha)^{-s}}{\Gamma(1-s)}}& \mbox{if}\quad x\in (\alpha,1] \\
\end{array}\right.
\end{equation}
and, for every $s\in(0,1)$, we have that
$$\int_\alpha^1 \frac{ (x-\alpha)^{-s}}{\Gamma(1-s)} = \left[\frac{(x-\alpha)^{1-s}}{\Gamma(2-s)}\right]_\alpha^1 =\frac{ (1-\alpha)^{1-s}}{\Gamma(2-s)}$$
which implies that
\begin{equation}\label{stric-cv}\|D^s_{0+}u \|_{L^1([0,1])} \overset{s\rightarrow 1}{\longrightarrow} 1 = |Du|([0,1])
\end{equation}
where $|Du|([0,1])$ denotes the total variation of $u$ on $[0,1]$.
\end{exa}
Next result is useful to prove Theorem \ref{WsSBV}:
  \begin{lemme}\label{w11-serie} Let $\{f_k\}\subset W^{1,1}(a,b)$ a sequence of non-negative functions with  $f_k(a)=0$ and
with non-negative derivative. We suppose also that
\begin{equation}\label{cond-deriv-serie}
\sum_k f_k\,,\; \sum_k f_k' \in \lone\,.
\end{equation}
Then
\begin{equation}\label{deriv-serie}
\left(\sum_k f_k(x)\right)' \,= \, \sum_k f_k'(x)\quad \quad  \mbox{ a.e.}\; x\in [a,b]\,.
\end{equation}
 \end{lemme}

 \noindent\textit{Proof.}
 The result follows from the monotone convergence theorem and the hypothesis on $f_k$. In fact, for every $x\in [a,b]$ we get
$$\sum_k f_k(x)= \sum_k  \int_a^x f_k'(t) \,dt =   \int_a^x \sum_k  f_k'(t) \,dt\,,$$
which proves that  $\sum_k f_k \in W^{1,1}(a,b)$  and \eqref{deriv-serie} follows.
 \proofend

\begin{thm} \label{WsSBV} For every $s \in (0,1) $,it holds that $SBV(a,b)\subset \wsp$ and
$$D^s_{a+} u  \overset{\mathcal{M}(a,b)}{\rightharpoonup} u'\mathcal{L}^1 + u(a^+) \delta_a +  \underset{x_k\in J_u}{\sum}p_k \delta_{x_k}  \,\quad  \mbox{as}\,\quad  s\rightarrow 1\, .$$
Moreover, if $ u(a^+)= 0 $ then
$$
\|D^s_{a+} u\|_{L^1(a,b)} \rightarrow \|u\|_{SBV(a,b)}\,\quad  \mbox{as}\,\quad  s\rightarrow 1\,,$$
where $u(a^+)$ denotes the right limit of $u$ at zero and  $J_u \subset (a,b) $ is the jump set of $u$. \end{thm}
\noindent\textit{Proof.}  With a simple change of variables, we can assume that $[a,b]=[0,1]$. This will  make the proof easier to read.  Every $BV$-function can be written as the difference of two increasing functions. Then, in the following we prove the result for  a $SBV$-increasing function.  Every $SBV$-function $u$ can be written as
$$u(x)=u(0^+) + \int_0^xu'(t)dt + \underset{x_k\in J_u}{\sum} p_k\chi_{[x_k,1]}(x) \quad \quad x\in[0,1]\,$$
where $u(0^+)$ denotes the right limit of $u$ at zero, $J_u$ denotes the (at most countable) set of jump points of $u$ and $p_k=u^+(x_k)-u^-(x_k)$ denotes the positive jump of $u$ at $x_k$. In particular
$$\|u\|_{SBV([0,1])}= |u(0^+)|+\int_0^1|u'(t)|dt + \underset{x_k\in J_u}{\sum} p_k\,.$$
In particular, $u$ can be written as the sum of two functions
$$
\begin{array}{l}
u(x)=f(x)+g(x)\,,\\
\displaystyle{f(x)= u(0^+) + \int_0^xu'(t)dt}\,,\\
\displaystyle{g(x)= \underset{x_k\in J_u}{\sum} p_k\chi_{[x_k,1]}(x) }\,.\\
\end{array}
$$
Now $f$  belongs to $W^{1,1}$ and using Theorem \ref{w11}, its fractional derivative is given by
$$f_s(x)=\frac{u(0^+)x^{-s}}{\Gamma(1-s)}+  I^{1-s}_{0+} [u'](x)\,.$$
So we  have $f_s\in L^1([0,1])$
and
\begin{equation}\label{conv-fs}
\begin{array}{lll}
\mbox{if} \; u(0^+)= 0 &
f_s  \overset{L^1(0,1)}{\longrightarrow} u' & \mbox{as}\,\quad  s\rightarrow 1\,,\\
\mbox{otherwise}&
f_s  \overset{\mathcal{M}(0,1)}{\rightharpoonup} u'\mathcal{L}^1  + u(0^+) \delta_0
 & \mbox{as}\,\quad  s\rightarrow 1\,.\\
\end{array}
\end{equation}
Concerning $g$, by the monotone convergence theorem, we have
$$I_{0+}^{1-s}g=  \underset{x_k\in J_u}{\sum} I_{0+}^{1-s}[ p_k\chi_{[x_k,1]}](x)\,.$$
Now, because of \eqref{charactder} and \eqref{charactder-bis}, the previous series verify the hypothesis of Lemma \ref{w11-serie}.  Condition \eqref{cond-deriv-serie} can be verified by observing that
$$\underset{x_k\in J_u}{\sum} I_{0+}^{1-s} [ p_k\chi_{[x_k,1]}]  \leq \underset{x_k\in J_u}{\sum} \frac{p_k}{\Gamma(2-s)}$$
and $\underset{x_k\in J_u}{\sum}  D_{0+}^s p_k\chi_{[x_k,1]}$ is well defined for every $x\notin J_u$ and
$$\int_0^1 \underset{x_k\in J_u}{\sum}  D_{0+}^s p_k\chi_{[x_k,1]}(t) \,dt =   \underset{x_k\in J_u}{\sum} I_{0+}^{1-s} [ p_k\chi_{[x_k,1]}](1) \leq \underset{x_k\in J_u}{\sum} \frac{p_k}{\Gamma(2-s)}\,.$$
With the previous proposition, the fractional derivative of $g$ is given by
$$g_s(x) =    (I_{0+}^{1-s}[g])' = \underset{x_k\in J_u}{\sum} \frac{ (x-x_k)^{-s}}{\Gamma(1-s)}  \chi_{[x_k,1]}\,.$$
Moreover, by  applying the monotone convergence theorem, we get
$$\|g_s\|_{\lone}= \underset{x_k\in J_u}{\sum}p_k\int_{x_k}^1 \frac{(t-x_k)^{-s}}{\Gamma(1-s)}dt = \underset{x_k\in J_u}{\sum}p_k\frac{(1-x_k)^{1-s}}{\Gamma(2-s)}
$$
and, similarly to \eqref{stric-cv} and using the fact that the series in the right-side is normally convergent, we have
\begin{equation}\label{conv-gs}
\|g_s\|_{\lone}\overset{s\rightarrow 1}{\longrightarrow}  \underset{x_k\in J_u}{\sum}p_k \,.
\end{equation}
The *-weak convergence of $g_s$ towards $ \underset{x_k\in J_u}{\sum}p_k \delta_{x_k}$ is obtained by the same arguments by writing every function test as the sum of its positive and negative part.
Finally, as the fractional derivative is a linear operator, we get
$$D^s_{0+}u(x) =  f_s(x)+g_s(x) \quad \forall \, x\in [0,1]\,.$$

Then, $D^s_{0+}u$ is defined at every point and belongs to $\lone$, which implies that $u\in W^{s,1}_{RL, 0+}([0,1])$ for every $s\in (0,1)$.
The result ensues from \eqref{conv-fs} and \eqref{conv-gs}.
\proofend

\begin{rem} We may note that if we extend  the functions $u$ with support in $ [a,b]$  by $0$  below $a$ and denote them similarly, we may consider $(-\infty, b]$ instead of $[a,b]$. Then the appearance of the Dirac measure at $a$ is consistent with the distributional derivative  of  $u$ on $(-\infty, b)$ since there is a jump at $a$.
\end{rem}
The subsequent remarks point out that the fractional Sobolev spaces are larger than $SBV$ and give some relationship between $BV$ and $W^{s,1}$.

\begin{rem}[{\bf Cantor-Vitali function}] The Cantor-Vitali function  is an example of increasing continuous function on $[0,1]$ whose standard derivative is defined and null at  a.e. point.  It is well known that such a function is of bounded variation but is not a $SBV$-function.  Precisely such a function is H\"older-continuous with exponent $\alpha = \ln 2/\ln 3$ (i.e., the Hausdorff dimension of the Cantor set). Then, with Proposition \ref{exist-der-hold}   the Riemann-Liouville derivative is well defined at every point for every $s\in (0,\alpha)$. In fact, it is  surely possible to prove the result for any $s \in (0,1)$  using an adapted Cantor-like  function whose H\"older exponent  would be $s$. Moreover the Riemann-Liouville derivative of order $s$ belongs to
$\mathcal{C}^{0,\alpha-s}([0,1])\,$, and we get that the Cantor-Vitali function belongs to $W^{s,1}_{RL, 0+} ([0,1])$ for every $s\in (0,\alpha)$. This proves in particular that in general
$$(BV\setminus SBV) \cap W^{s,1}_{RL, 0+} \neq \emptyset\,,$$
since the function we exhibit belongs to  $BV$ and  $W^{s,1}_{RL, 0+} $ and not  to $SBV$. \end{rem}

\begin{rem}[{\bf A continuous but non H\"older-continuous function in RL spaces}]
  Set $u(x)=\big(\ln(x/2)\big)^{-1}$\,if $ x\in(0,1)$\, and $u(0)=0\,.$
  This $u$ provides an example of monotone, continuous function 
    which is not $\alpha$ H\"older-continuous for any $\alpha\in(0,1)$,
    but $u$ belongs to $\underset{s\in(0,1)}{\bigcap} W^{s,1}$\,.
\end{rem}

\begin{rem}[{\bf Relationship between $BV$ and  $W^{s,1}_{RL, 0+} $  }] In \cite{samko_exmp} the authors investigate the relationship between usual a.e. differentiation and the fractional Riemann-Liouville derivative definition. Several interesting examples are given.
One of them  is given by the Weierstrass function defined as
$$W(x) =\sum_{n=0}^\infty q^{-n}(e^{iq^{n}x}-e^{iq^{n}a}) \quad \quad  x\in [a,b] $$
where $q>1$. It is proved that  $W$  has continuous and
bounded fractional Riemann-Liouville derivatives of all orders $s < 1.$  However, since  $W$ is nowhere differentiable it cannot be of bounded variation.

This implies that Riemann-Liouville  fractional Sobolev spaces are  not contained in $\BV$. Then we can state
$$SBV\subset \underset{s\in(0,1)}{\bigcap} W^{s,1}\;,\quad \quad \underset{s\in(0,1)}{\bigcap} W^{s,1}\setminus BV \neq \emptyset\,.$$
The question  to be addressed now is the relation with
$BV $ and $W^{s,1}_{RL, 0+} $. Indeed, we have either $BV \subset W^{s,1}_{RL, 0+} $ or they are completely different spaces  whose intersection contains $SBV$ (the case $W^{s,1}_{RL, 0+}\subset BV$ is excluded).

\end{rem}

 \section{Conclusion}

In this paper we try to make connections between the two main definitions of fractional derivatives : the \textit{local} (pointwise) one whose typical representation is the RL derivative and the  \textit{global} one which is typically the Gagliardo one. In view of a more precise description of the derivative of order $s \in (0,1)$ with  respect to $BV$ functions we have also proved preliminary results to compare  $SBV(a,b) $ and $ \wsoneRL$.
Open problems are numerous. In particular, it remains to strongly connect the Riemann-Liouville theory with the Gagliardo one. This would allow to perform comparison between $ \wsoneRL$  and  Besov-spaces for example. In addition, we have to understand precisely how $ \wsoneRL$   behaves with respect  to $BV(a,b)$ to get some  information  about the $BV\backslash W^{1,1}$  functions.  In particular, we proved that $SBV(a,b)\subset \wsoneRL$  and exhibit a function  in  $\wsoneRL$, that does not  belong to $ BV(a,b)$; however we still don't know if $ BV(a,b) \subsetneq \wsoneRL  $.
  In addition, it remains to prove density results and continuity/compactness results in view of variational models involving the RL derivative.
  
  Another  important issue is also  to address engineering and/or imaging problems  in a rigorous mathematical framework. From that point of view, the fractional derivative concept  widely used in  engineering is the Caputo  one : 
\begin{defn}[{\bf Caputo fractional derivative}]\cite{almeida,caputo}
Let $u\in L^1(a,b)$ and $n-1\leq s < n$ (where $n$ integer).\\
The left Caputo fractional derivatives of $u$ at $x\in \set$ is defined by
\begin{equation}\label{caputo}
 ^C\! D_{a+}^s u(x)\, =\, I_{a+}^{n-s} \left[\frac{d^n}{d x^n}  u\right] (x) =   \frac{1}{\Gamma(n-s)} \int_a^x   \frac{u^{(n)}(t)}{(x-t)^{s-n+1}} d t\,
 \end{equation}
\vskip-0.1cm \noindent when the right-hand side is defined.  
The right Caputo fractional derivatives of $u$ at $x\in \set$ is defined  in a similar way. 
\end{defn}
The main advantage of Caputo derivatives, which makes them the preferred ones in many engineering applications,
is the fact that the initial conditions for fractional differential equations with Caputo derivatives 
are expressed by integer derivatives at time 0, say quantities with a straightforward physical interpretation.
This relies on the Laplace transform of Caputo fractional derivative, formally identical to the classical formula for integer derivatives (in contrast to the formula for RL fractional derivatives):
\begin{equation*}\label{LaplaceCaputo}
  \mathcal{L}\left\{\, ^C\! D_{a+}^\alpha u \,\right\}(s)\ =\ s^\alpha \,\mathcal{L}\{u\}\,-\, \sum_{j=0}^{n-1}\,s^{\alpha-j-1}u^{j}(0) \qquad n-1<\alpha\leq n \, .
\end{equation*}
The connection between Riemann-Liouville and Caputo derivatives, when they both exist, is given by the relationship (see\cite{almeida} for example):\vskip-0.3cm
\begin{equation*}\label{left RL-C}
^C\! D_{a+}^\alpha u(x)\, =\, D_{a+}^\alpha u(x)\,-\, \sum_{j=0}^{n-1}\,\frac {u^{j}(a)}{\Gamma(j-\alpha+1)}\,(x-a)^{j-\alpha}\,.
  \end{equation*}
In particular:
\begin{equation*}\label{left RL=C}
^C\! D_{a+}^\alpha u(x)\, =\, D_{a+}^\alpha u(x)\qquad \hbox{if } u(a)=u'(a)=\dots=u^{n-1}(a)=0\,.
  \end{equation*}
Therefore all the results and comments stated in the present paper concerning Riemann-Liouville derivatives 
can be easily transferred to Caputo derivatives. This will be precisely addressed in a forthcoming work. 
%

\bibliographystyle{plain}  
\bibliography{bibliofrac}	

 \bigskip \smallskip

 \it

 \noindent
$^1$  Laboratoire MAPMO, CNRS, UMR 7349,\\
 F\'ed\'eration Denis Poisson, FR 2964,\\
  Universit\'{e} d'Orl\'{e}ans, 
 B.P. 6759,\\
 45067 Orl\'{e}ans cedex 2,
France \\[4pt]
e-mail: {maitine.bergounioux@univ-orleans.fr}  
\\[12pt]
$^2$  Universit\`a del Salento, \\Dipartimento di Matematica e Fisica ``Ennio De Giorgi'',\\
I 73100 Lecce, Italy\\
[4pt]
e-mail: antonio.leaci@unisalento.it
\\[12pt]
$^3$ Institut Pasteur, Laboratoire d'Analyse d'Images Biologiques, \\ CNRS, \ UMR \ 3691, \\  Paris, \ France,  \\[4pt]
e-mail:giacomo.nardi@pasteur.fr 
\\[12pt]
$^4$ Politecnico di Milano, Dipartimento di Matematica,\\
Piazza ``Leonardo da Vinci'', 32,\\
I 20133 Milano, Italy \\
[4pt]
e-mail:{franco.tomarelli@polimi.it}

\end{document}